\documentclass[aos,MSNbibl,seceqn,citesort,dvips]{arximspdf}
\doi{10.1214/11-AOS964}
\volume{40}
\issue{1}
\pubyear{2012}
\firstpage{412}
\lastpage{435}

\makeatletter

\newcommand{\bigtimes}{\mathop{\mbox{\fontsize{17}{17}\selectfont{$\!\times$}}}}

\newcommand{\cal}{\mathcal}

\newtheorem{theorem}{Theorem}[section]
\newtheorem{prp}{Proposition}[section]
\newtheorem{lem}{Lemma}[section]
\newtheorem{cor}{Corollary}[section]

\newproclaim{Remark}{Remark}

\newcommand{\Var}{\operatorname{Var}}
\newcommand{\Cov}{\operatorname{Cov}}
\renewcommand{\P}{\mathbb{P}}

\newcommand{\eps}{\varepsilon}

\newcommand{\cN}{{\cal N}}
\newcommand{\cC}{{\cal C}}

\newcommand{\bbR}{\mathbb{R}}

\newcommand{\field}[1]{\mathbb{#1}}
\newcommand{\R}{\field{R}}
\newcommand{\E}{\field{E}}
\newcommand{\EXP}{\E}
\newcommand{\PROB}{\field{P}}
\newcommand{\bA}{\mathbf{A}}
\newcommand{\bB}{\mathbf{B}}
\newcommand{\bI}{\mathbf{I}}
\newcommand{\bM}{\mathbf{M}}

\makeatother

\begin{document}
\begin{frontmatter}

\title{Detection of correlations}
\runtitle{Detection of correlations}

\begin{aug}
\author[A]{\fnms{Ery} \snm{Arias-Castro}\thanksref{t1}\ead[label=e1]{eariasca@math.ucsd.edu}},
\author[B]{\fnms{S{\'e}bastien} \snm{Bubeck}\corref{}\ead[label=e2]{sbubeck@princeton.edu}}
\and
\author[C]{\fnms{G{\'a}bor} \snm{Lugosi}\thanksref{t2}\ead[label=e3]{gabor.lugosi@upf.edu}}
\runauthor{E. Arias-Castro, S. Bubeck and G. Lugosi}
\affiliation{University of California, San Diego,
Princeton University, and ICREA and Pompeu Fabra University}
\address[A]{E. Arias-Castro\\
Department of Mathematics\\
University of California, San Diego\\
San Diego, California 92093\\
USA\\
\printead{e1}}
\address[B]{S. Bubeck\\
Department of Operations Research\\
\quad and Financial Engineering\\
Princeton University\\
Princeton, New Jersey 08542\\
USA\\
\printead{e2}}
\address[C]{G. Lugosi\\
ICREA\\
and\\
Department of Economics\\
Pompeu Fabra University\\
Barcelona\\
Spain\\
\printead{e3}} 
\end{aug}

\thankstext{t1}{Supported by ONR Grant N00014-09-1-0258.}

\thankstext{t2}{Supported by the Spanish Ministry of Science
and Technology Grant MTM2009-09063
and \mbox{PASCAL2} Network of Excellence under EC Grant 216886.}

\received{\smonth{9} \syear{2011}}
\revised{\smonth{12} \syear{2011}}

%
\begin{abstract}
We consider the hypothesis testing problem of deciding whether an
observed high-dimensional vector has independent normal components or,
alternatively, if it has a~small subset of correlated components.
The correlated components may have a~certain combinatorial structure
known to the statistician. We establish upper and lower bounds for
the worst-case
(minimax) risk in terms of the size of the correlated subset, the level
of correlation, and the structure of the class of possibly correlated
sets. We show that some simple tests have near-optimal performance in
many cases, while the generalized likelihood ratio test is suboptimal
in some important cases.
\end{abstract}

%

\end{frontmatter}

\section{Introduction}\label{intro}

In this paper we consider the following statistical problem:
upon observing a~high-dimensional vector, one is interested
in detecting the presence of a~sparse, possibly structured,
correlated subset of components of the vector. Such problems emerge naturally
in numerous scenarios.
The setting is closely related to Gaussian signal detection in Gaussian
white noise, on which there is an extensive literature surveyed in
\cite{MR1658799}. In image processing, textures are modeled via Markov
random fields~\cite{4767341}, so that detecting a~textured object
hidden in Gaussian white noise amounts to finding an area in the image
where the pixel values are correlated. Similar situations arise in
remote sensing based on a~variety of hardware. A~related task is the
detection of space--time correlations in multivariate time series, with
potential applications to finance~\cite{cboe}.

\subsection{Setting and notation}
\label{secsetting}

We investigate the possibilities and limitations in
problems of detecting correlations in a~Gaussian framework.\vadjust{\goodbreak}
We may formulate this as a~general hypothesis testing problem as follows.
An $n$-dimensional
Gaussian vector $X=(X_1,\ldots,X_n)$ is observed. Under the null hypothesis
$H_0$, the vector
$X$ is standard normal, that is, with zero mean vector and identity
covariance matrix. To describe the alternative hypothesis~$H_1$,
let $\cC$ be a~class of subsets of $\{1,\ldots, n\}$, each of size $k$,
indexing the possible ``contaminated'' components.
One wishes to test whether there exists an $S \in\cC$ such that
\[
\Cov(X_i, X_j) = \cases{
1, &\quad $i = j$, \cr
\rho, &\quad $i \neq j$, with $i, j \in S$, \cr
0, &\quad otherwise,}
\]
where $\rho>0$ is a~given parameter.
Equivalently, if $X = (X_1,\ldots, X_n)$ denotes the vector of
observations, then
\[
H_0\dvtx X \sim\cN(0,\bI)\quad\mbox{vs.}\quad
H_1\dvtx X \sim\cN(0,\bA_S)\qquad\mbox{for some $S \in\cC$,}
\]
where
$\bI$ denotes the $n\times n$ identity matrix and
%
\begin{equation} \label{model}
(\bA_S)_{i,j} = \cases{
1, &\quad $i = j$, \cr
\rho, &\quad $i \neq j$, with $i, j \in S$, \cr
0, &\quad otherwise.}
\end{equation}
We write $\PROB_0$ for the probability under $H_0$ (i.e., the standard
normal measure in $\R^n$) and, for each $S\subset\cC$,
$\PROB_S$ for the measure of $\cN(0,\bA_S)$.

The goal of this paper is to understand for what values of the
parameters $(n,k,\rho)$ reliable testing is possible. This, of course, depends
crucially on the size and structure of the subset class $\cC$. We
consider the following two
prototypical
classes:
\begin{itemize}
\item\textit{$k$-intervals.}
In this example, we consider the class of all intervals of size $k$ of
the form $\{i,\ldots, i + k -1\}$ modulo $n$---for
aesthetic reasons. (We call such an interval a~\textit{$k$-interval}.)
This class is the flagship of \textit{parametric} classes, typical of
the class of objects of interest in signal processing.
%
\item\textit{$k$-sets.}
In this example, we consider the class of all sets of size $k$, that
is, of the form $\{i_1,\ldots, i_k\}$ where the indices are all
distinct in $\{1,\ldots, n\}$. (We call such a~set a~\textit{$k$-set}.) This class is the flagship of \textit{nonparametric}
classes, and may arise in multiple comparison situations.
\end{itemize}
Our theory, however, applies more generally to other classes, such as:
\begin{itemize}
\item\textit{$k$-hypercubes.} In this example, the variables are
indexed by the $d$-dimen\-sional lattice, that is, $X = (X_{i}\dvtx i \in
\{1,\ldots, m\}^d)$, so that the sample size is $n = m^d$, and we
consider the class of all hyper-rectangles of the form
$\bigtimes_{s=1}^d \{i_s,\ldots, i_s + k_s -1\}$---each interval modulo
$m$---of fixed size $\prod_{s=1}^d k_s = k$. This class is the
simplest model for objects to be detected in images (mostly $d=2,3$
in applications).
\item\textit{Perfect matchings.}
Suppose $n$ is a~perfect square with $k^2=n$. The
components of the observed vector $X$ correspond to edges of the
complete bipartite graph on $2k$ vertices and each set in $\cC$ corresponds
to the edges of a~perfect matching. Thus, $|\cC|=k!$.
In this example $\cC$ has a~nontrivial combinatorial structure.
%
\item\textit{Spanning trees.}
In another example, $n={k+1\choose2}$ and the components of $X$ correspond
to the edges of a~complete graph $K_{k+1}$ on $k+1$ vertices and every
element of $\cC$ is a~spanning tree of $K_{k+1}$.
\end{itemize}

As usual, a~\textit{test} is a~binary-valued function $f\dvtx\R^n \to\{
0,1\}$. If
$f(X)=0$, then the test accepts the null hypothesis $H_0$;
otherwise $H_0$ is rejected by $f$.
We measure the performance of a~test based on its \textit{worst-case
risk} over the class of interest $\cC$, formally defined by
\[
R^{\max}(f) = \PROB_0\{f(X)=1\}
+ \max_{S\in\cC} \PROB_S\{f(X)=0\}.
\]
%
We will derive upper and lower bounds on the \textit{minimax risk}
\[
R_*^{\max} := \inf_{f} R^{\max}(f).
\]
%
A~standard way of obtaining lower bounds for the minimax risk
is by putting a~prior on the
class $\cC$ and obtaining a~lower bound on the corresponding \textit
{Bayesian risk}, which never exceeds the worst-case risk. Because
this is true for any prior, the idea is to find one that is hardest
(often called \textit{least favorable}).
Most
classes we consider here
are invariant under some group action: $k$-intervals are invariant
under translation and $k$-sets are invariant under permutation.
Invariance considerations (\cite{TSH}, Section 8.4) lead us to
considering the
uniform prior on $\cC$, giving rise to the following \textit{average risk}:
\[
R(f) = \PROB_0\{f(X)=1\}
+ \PROB_1\{f(X)=0\},
\]
where
\[
\PROB_1\{f(X)=0\} := \frac{1}{N}\sum_{S\in\cC} \PROB_S\{f(X)=0\},
\]
and $N := |\cC|$ is the cardinality of $\cC$.
%
The advantage of considering the average risk over the worst-case
risk is that we know an optimal test for the former, which, by the
Neyman--Pearson fundamental lemma, is the likelihood ratio test,
denoted $f^*$.
%
Introducing
%
\begin{equation} \label{Z}
Z_S = \exp\bigl(\tfrac{1}{2} X^T (\bI- \bA_S^{-1}) X \bigr)
\end{equation}
for all $S \in\cC$,
the likelihood ratio between $H_0$ and $H_1$
may be written as
%
\begin{equation} \label{L}
L(X) = \frac{1}{N} \sum_{S \in\cC} \frac{Z_S}{\E_0 Z_S},
\end{equation}
and the optimal test becomes
\[
f^*(x) = 0  \quad\mbox{if and only if}\quad   L(x) \le1.
\]
Note that $\E_0 Z_S = \sqrt{\det(\bA_S)}$.
The
(average)
risk $R^*=R(f^*)$ of the optimal test is called the
\textit{Bayes risk} and it satisfies
\[
R^* = 1 - \frac{1}{2} \E_0 |L(X) - 1| = 1 - \frac{1}{2} \E_0
\biggl| \frac{1}{N} \sum_{S \in\cC} \frac{Z_S}{\E_0 Z_S} - 1\biggr|.
\]
%
Note that, with the only exception of the case of spanning trees,
in all examples mentioned above,
the minimax and Bayes
risks coincide, that is, $R^* = R_*^{\max}$. This is again due to
invariance (\cite{TSH}, Section 8.4).
(The class of spanning trees is not sufficiently symmetric for
this equality to hold. However, as we will see below, even in this
case, $R^*$ and $R_*^{\max}$ are of the same order of magnitude.)


We focus on the case when $n$ is large and formulate some
of the results in an asymptotic language with $n\to\infty$ though
in all cases explicit nonasymptotic inequalities are available.
Of course, such asymptotic statements only make sense if we define
a~sequence of integers $k=k_n$ and classes $\cC=\cC_n$. This
dependency in $n$ will be left implicit.
In this asymptotic setting, we say that \textit{reliable} detection is
possible (resp., impossible) if
$R_*^{\max}\to0$
(resp., $\to1$) as $n \to\infty$.

%
\begin{Remark*}[\,(Covariance structure)]
In this paper we assume that, under the alternative hypothesis,
the correlation between any two variables in the ``contaminated'' set
is the same. While this model has a~natural interpretation (see
Lem\-ma~\ref{lemrepresent} below), it is clearly a~restrictive assumption.
This simplification is \mbox{convenient} in understanding the fundamental
limits of detection (i.e., in obtaining lower bounds on the risk).
At the same time, the tests we exhibit also match these lower bounds
under more general correlation structures, such as
%
\begin{equation} \label{model-general}
(\bA_S)_{i,j} \cases{
= 1, &\quad $i = j$, \cr
\geq\rho, &\quad $i \neq j$, with $i, j \in S$,
\cr
= 0, &\quad otherwise.}
\end{equation}
That said, dealing with more general correlation structures remains an
interesting and important challenge, relevant in the detection of
textured objects in textured background, for example.
\end{Remark*}

\subsection{Relation to previous work}

The vast majority of the literature on detection is concerned with the
detection of a~signal in additive (often Gaussian) noise, which would
correspond here to an alternative where $X_i \sim\cN(\mu, 1)$ for $i
\in S$, where $\mu> 0$ is the (per-coordinate) signal amplitude. We
call this the \textit{detection-of-means} setting. The literature on
this problem is quite comprehensive. Indeed, the detection of
$k$-intervals and $k$-hypercubes is treated extensively in a~number of
papers; see, for example,~\cite{MGD,morel,perone,boutsikas,cluster}. A~more general framework that includes the detection of perfect
matchings and spanning trees is investigated in~\cite{combin}, and the
detection of $k$-sets is studied
in~\cite{JinPhD,Ingster99,baraud,dj04,hj09}.
In the literature on detection of parametric objects, the phrase
``correlation detection'' usually refers to the method of \textit{matched
filters}, which consists of correlating the observed signal with
signals of interest. This is not the problem we are interested in
here.
While the problem of \textit{detection-of-correlations} considered here
is mathematically more challenging than the detection-of-means setting,
there is a~close relationship between the two.
The connection is established by the representation theorem
of~\cite{MR0145564}---stated here for the case Gaussian random
variables.
%
\begin{lem}[(\cite{MR0145564})] \label{lemrepresent}
Let $X_1,\ldots, X_k$ be standard normal with $\Cov(X_i, X_j) = \rho
$ for $i \neq j$. Then there are i.i.d. standard normal random
variables, denoted $U, U_1,\ldots, U_k$, such that $X_i = \sqrt{\rho}
U + \sqrt{1-\rho}   U_i$ for all $i$.
\end{lem}

Thus, given $U$, the problem becomes that of detecting a~subset of
variables with nonzero mean (equal to $\sqrt{\rho}   U$) and with a~variance equal to $1-\rho$ (instead of 1). This
simple observation will be very useful to us later on. When $U$ is random,
the setting is similar to that of detecting a~Gaussian process (here
equal to $\sqrt{\rho}   U$ for $i \in S$, and equal to 0 otherwise)
in additive Gaussian noise. However, the typical setting assumes that
the Gaussian process affects all parts of the
signal~\cite{MR1658799}. In our setting, the signal (the subset of
correlated variables) will be sparse.
Since we only have one instance of the signal $X$, the problem cannot
be considered from the perspective of either multivariate statistics
or multivariate time series. If indeed we had multiple copies of $X$,
we could draw inspiration from the literature on the estimation of
sparse correlation matrices~\cite{MR2485008,MR2676885}, from the
literature on multivariate time series~\cite{MR2722660}, or
on other approaches~\cite{DeGyLuUd11};
but this is
not the case as we only observe $X$.
Closer in spirit to our goal of detecting correlations in a~single
vector of observation is the paper of~\cite{MR2597269}, which aims at
testing whether a~Gaussian random field is i.i.d. or has some Markov
dependency structure. Their setting models communication networks and
is not directly related to ours.

It transpires, therefore, that $\rho$ in the detection-of-correlations
setting plays a~role analogous to $\mu^2$ in the detection-of-means
setting. While this is true to a~certain extent, the picture is quite
a~bit more subtle. The detection-of-means problem for parametric
classes such as $k$-intervals is well understood. In such cases,
$\mu^2$~needs to be of order at least $(1/k)
\log(n/k)$ for reliable detection of $k$-intervals to be possible.
This remains true in the detection-of-correlations setting, and the
\textit{generalized likelihood ratio test} (\textit{GLRT}) is near-optimal,
just as in the detection-of-means problem; see, for example,~\cite{MGD}.

Our inspiration for considering $k$-sets comes from the line of
research on the detection of sparse Gaussian mixtures. Very precise
results are known on $(n,k,\mu)$ that make detection
possible~\cite{JinPhD,Ingster99,baraud} and optimal tests have been
developed, such as the ``higher criticism''~\cite{dj04,hj09}. In fact,
the recent paper~\cite{cai-jeng-jin} deals with heteroscedastic
instances of the detection-of-means problem where the variance of the
anomalous variables may be different from 1.
For example, it is known that, when $n = O(k^2)$ [resp., $k^2 = o(n)$],
$\mu^2$ needs to be of order at least $n/k^2$ [resp., $\log(n)$] for
reliable detection of $k$-sets to be possible, and the test based on
$\sum_i X_i$ (resp., $\max_i X_i$) is near-optimal. Though more
precise results are available when $k^2 = o(n)$, these cannot
be translated immediately to our case via the
representation theorem of Lemma~\ref{lemrepresent}.
As a~bonus, we show that the GLRT is clearly suboptimal in some
regimes---see Theorem~\ref{thmglrt-bad}. Note that in the
detection-of-means problem it is not known whether the GLRT has any power.

\subsection{Contribution and content of the paper}

This paper contains a~collection of positive and negative results about
the detection-of-correlation problem described above.
In Section~\ref{seclower} we derive lower bounds for the Bayes risk.
The usual route of bounding the
variance of the likelihood ratio, that is very successful in the
detection-of-means problem, leads essentially nowhere in our case.
Instead, we develop a~new approach based on Lemma~\ref{lemrepresent}.
We establish a~general lower bound for the Bayes
risk in terms of the moment generating function of the
size of
the overlap of two randomly chosen elements of the class $\cC$.
This quantity also plays a~crucial role
in the detection-of-means setting and we are able to use
inequalities worked out in the literature in various examples.
In Section~\ref{secupper} we study the performance of some simple and
natural tests such as the squared-sum test---based on $(\sum_i
X_i)^2$, the generalized likelihood ratio test (GLRT) and a~goodness-of-fit (GOF) test, as well as some variants.
We show that, in the case of parametric classes such as $k$-intervals
and $k$-hypercubes, the GLRT is essentially optimal. The squared-sum
test is shown to be essentially optimal in the case of $k$-sets when
$k^2/n$ is large, while the GLRT is clearly suboptimal in this
regime. This is an interesting example where the GLRT fails
miserably. When $k^2/n$ is small, detection is only possible
when $\rho$ is very close to $1$. We show that a~simple GOF test is
near-optimal in this case.
The analysis of tests such as the squared-sum test and the GLRT
involves handling quadratic forms in $X$. This is technically
more challenging than the analogous problem for
the detection-of-means setting in which only linear functions of $X$
appear (which are normal random variables).

\section{Lower bounds}
\label{seclower}

In this section we investigate lower bounds on the risk, which are
sometimes called information bounds. First
we consider the special case when $\cC$ contains only one element as
this example will serve as a~benchmark for other examples.
Then we consider the standard method based on bounding the variance of the
likelihood ratio under the null hypothesis, and show that it leads nowhere.
We then develop a~new bound
based on Lemma~\ref{lemrepresent} that
has powerful implications, leading to fairly sharp bounds in a~number of examples.

\subsection{The case $N=1$}
\label{secsimple}

As a~warm-up, and to gain insight into the problem,
consider first the simplest case where $\cC$
contains just one set, say $S=\{1,\ldots,k\}$.
In this case, the alternative hypothesis is simple
and the likelihood ratio (Neyman--Pearson)
test may be expressed by
\[
f^*(X) =0  \quad\mbox{if and only if}\quad
X^T (\bI- \bA_S^{-1}) X \le\log\det(\bA_S).
\]
%
This follows by the fact that $\EXP Z_S= \sqrt{\det(\bA_S)}$
which is easy to check by straightforward calculation.

The next simple lemma helps understand the behavior of the Bayes risk.
%
\begin{lem}
\label{lemqf}
Under $\PROB_0$, $X^T (\bI- \bA_S^{-1}) X$ is distributed as
\[
- \frac{\rho}{1-\rho} \chi^2_{k-1} + \frac{\rho(k-1)}{1 + \rho
(k-1)} \chi^2_1,
\]
and under the alternative $\PROB_S$, it has the same distribution as
\[
-\rho\chi^2_{k-1} + \rho(k-1) \chi^2_{1},
\]
where $\chi_1^2$ and $\chi_{k-1}^2$ denote independent $\chi^2$
random variables with degrees of freedom $1$ and $k-1$, respectively.
\end{lem}
\begin{pf}
If $Y = (Y_1,\ldots, Y_n)$ denotes a~standard normal vector, then
under $H_0$, the quadratic form $X^T (\bI- \bA_S^{-1}) X$
is distributed as $Y^T (\bI- \bA_S^{-1}) Y$,
and under
the alternative, it has the distribution of
$Y^T (\bA_S - \bI) Y$, since $X$ is distributed as $\bA_S^{1/2} Y$.

Now, observe that for any symmetric matrix $\bB$ with eigenvalues
$\lambda_1,\ldots,\lambda_n$, the quadratic form $Y^T \bB Y$ has distribution
%
\begin{equation} \label{sum-chi}
Y^T \bB Y \sim\sum_{i=1}^n \lambda_i Y_i^2.
\end{equation}
This follows simply by diagonalizing $\bB$ and using the rotational
invariance of the standard normal distribution.

The lemma follows from this simple representation and the
fact that
$\bA_S$ has eigenvalue $1-\rho$ with
multiplicity $k-1$, $1 + \rho(k-1)$ with multiplicity~$1$, and the
eigenvalue $1$ with multiplicity $n-k$.
\end{pf}

Now it is straightforward to analyze the Bayes risk. In particular,
we immediately have the following:
%
\begin{prp}
\label{prpsimple}
If $\cC$ is a~singleton, $\lim_{k\to\infty}R^*= 0$ if and only if
$\rho k \to\infty$. Similarly,
$\lim_{k\to\infty}R^*= 1$ if and only if
$\rho k \to0$.
\end{prp}
\begin{pf}
Suppose $\rho k \to\infty$.
It suffices to show that there exists a~threshold $\tau_k$ such that
$\PROB_0\{X^T (\bI- \bA_S^{-1}) X\ge\tau_k\} \to0$
and $\PROB_S\{X^T (\bI- \bA_S^{-1}) X< \tau_k\} \to0$.
We use Lemma~\ref{lemqf} and the fact that, by Chebyshev's inequality,
\[
\mathbf{P}\bigl\{|\chi_{k}^2 - k| > t_k \sqrt{k}\bigr\} \to0, \qquad
k \to \infty,
\]
for any sequence $t_k \to\infty$, and the fact that
\[
\mathbf{P}\{t_k^{-1} < \chi_1^2 < t_k\} \to1  \qquad\mbox{as $k \to
\infty$}.
\]
We choose $t_k = \log k$ and define $\tau_k := -\rho k + \rho t_k
\sqrt{k} + t_k$.
Then under the null,
\[
\PROB_0\{X^T (\bI- \bA_S^{-1}) X \ge\tau_k \} \to0,
\]
and under the alternative, setting $\eta_k := -\rho k - \rho t_k \sqrt
{k} + \rho k t_k^{-1}$,
\[
\PROB_S\{X^T (\bI- \bA_S^{-1}) X < \eta_k \} \to0.
\]
We then conclude with the fact that, for $k$ large enough, $\tau_k <
\eta_k$.

If $\rho k$ is bounded, the densities of the test statistic under
both hypotheses have a~significant overlap and the risk cannot
converge to $0$.

The proof of the second statement is similar.
\end{pf}

Clearly, the role of $n$ is immaterial in this specific example as
the optimal test ignores all components whose indices are not in $S = \{
1,\ldots, k\}$.

\subsection{The moment method}
\label{secmoment}

When the class $\cC$ contains more than one element, the likelihood
ratio with uniform prior on $\cC$ is given by~(\ref{L}). A~common
approach for deriving a~lower bound on the Bayes risk is via an upper
bound on the \textit{variance} of $L(X)$ under the null. Indeed,
by the Cauchy--Schwarz inequality,
\[
R^* = 1 - \frac{\E_0 |L(X)-1|}{2} \geq1 - \frac{\sqrt{\E_0
[L(X)^2] - 1}}{2}.
\]
Therefore, an upper bound on $\E_0 [L(X)^2] -1 = \Var_0(L(X))$ leads
to a~lower bound on $R^*$.

Let $\Lambda= \det(\bA_S) = (1-\rho)^{k-1} (1 + \rho(k-1))$,
which is independent of $S \in\cC$.
By Fubini's theorem, we have
\[
\E_0 L(X)^2 = \frac1\Lambda\frac1{N^2} \sum_{S, S' \in\cC} \E
_0(Z_S Z_{S'}),
\]
where $Z_S$ is defined in~(\ref{Z}). We focus on terms of the double
sum for which $S = S'$.\vadjust{\goodbreak}

The following result is a~straightforward consequence of the
representation~(\ref{sum-chi}) and the well-known expression for the
moment generating function of~$\chi_1^2$.
%
\begin{lem} \label{lemchi-eig}
Suppose $X$ is a~standard normal vector in $\bbR^n$ and $\bM$ is an
$n\times n$ symmetric matrix with eigenvalues strictly less than $1/2$. Then
\[
\E\exp(X^T \bM X) = \det(\bI- 2\bM)^{-1/2}.
\]
If $\bM$ has an eigenvalue exceeding $1/2$, then
$\E\exp(X^T \bM X) = +\infty$.
\end{lem}

Since $\bM:= \bI- \bA_S^{-1}$ has eigenvalue $-\rho/(1-\rho)$ with
multiplicity $k$, eigenvalue $\rho(k-1)/(1 + \rho(k-1))$ with
multiplicity $1$, and eigenvalue 0 with multiplicity $n - k$,
$\E_0 [Z_S^2] = \E_0 \exp(X^T \bM X) = +\infty$ unless
$\rho
(k-1) < 1$. The implications are rather insubstantial. It only shows
that, when $\rho(k-1) \leq1 -\eps$ with $\eps> 0$ fixed, the Bayes
risk does not tend to zero. As we shall see, this lower bound is
grossly suboptimal,
except in the case where $\cC$ is a~singleton (as in Section~\ref{secsimple})
or does not grow in size with $n$.

A~refinement of this method consists in bounding the first and second
\textit{truncated} moments of $L(X)$, again under the null hypothesis.
For example, this is the approach used in~\cite{Ingster99,cai-jeng-jin}
in the detection-of-means setting
for the case of $k$-sets to obtain sharp
bounds. Unfortunately, in our case this method only provides
a~useful bound when the class $\cC$ is not too large
(i.e., has size polynomial in $k$)
while it does not seem to lead anywhere in the case of $k$-sets.
The computations are quite involved and we do not
provide details here, as we were able to obtain a~more powerful
general bound that
applies to both $k$-intervals and $k$-sets. This is presented
in the next section.

\subsection{A~general lower bound}
\label{secgeneral}

In this section we derive a~general lower bound for the Bayes risk.
As in the
detection-of-means problem~\cite{maze,cluster,combin}, the relevant
measure of complexity is in terms of the moment generating function of
the size of the overlap of two randomly chosen elements of $\cC$. In
the detection-of-means setting, this is a~consequence of bounding the
variance of the likelihood ratio. We saw in Section~\ref{secmoment}
that this
method is useless here. Instead, we make a~connection between
the two problems using Lemma~\ref{lemrepresent}.
%
\begin{theorem}
\label{thmlower}
For any class $\cC$ and any $a~> 0$,
\[
R^* \geq\mathbf{P}\{|\cN(0,1)| \le a\} \bigl(1 - \tfrac12 \sqrt{\E
\exp
(\nu_a~Z ) - 1}\bigr),
\]
where $\nu_a~:= \rho a^2/(1+\rho)- \frac12 \log(1-\rho^2)$ and
$Z=|S \cap S'|$, with $S, S'$ drawn independently, uniformly at random
from $\cC$. In particular, taking $a~= 1$,
\[
R^* \geq0.6 - 0.3 \sqrt{\E\exp(\nu_1 Z ) - 1},
\]
where $\nu_1 = \nu(\rho) := \rho/(1+\rho)- \frac12 \log(1-\rho^2)$.
\end{theorem}
\begin{pf}
The starting point of the proof is
Lemma~\ref{lemrepresent},\setcounter{footnote}{2}\footnote{In
fact, we only need to assume that $X$ is as described in distribution.}
which enables us to represent the vector $X$ as
\[
X_i = \cases{
U_i, &\quad if $i \notin S$, \cr
\sqrt{\rho}  U + \sqrt{1-\rho}  U_i, &\quad if $i \in S$,}
\]
where $U,U_1,\ldots,U_n$ are independent standard normal random variables.

We consider now the alternative $H_1(u)$, defined as the alternative
$H_1$ given \mbox{$U=u$}.
Let $R(f)$, $L$, $f^*$ [resp., $R_u(f)$, $L_u$, $f_u^*$] be the risk of
a~test $f$, the likelihood ratio, and the optimal (likelihood ratio)
test, for $H_0$ versus $H_1$ [resp., $H_0$~versus $H_1(u)$]. For any $u
\in\bbR$, $R_u(f_u^*) \leq R_u(f^*)$, by the optimality of $f_u^*$
for $H_0$ versus $H_1(u)$. Therefore, conditioning on $U$,
\begin{eqnarray*}
R^* & = & R(f^*) \\
& = & \E_{U} R_U(f^*) \\
& \geq& \E_{U} R_U(f_U^*) \\
& = & 1 - \tfrac{1}{2} \E_{U} \E_0 |L_U(X) - 1|.
\end{eqnarray*}
[$\E_U$ is the expectation with respect to $U \sim\cN(0,1)$.]
Using the fact that
$\E_0 |L_u(X) - 1| \le2$ for all $u$, we have
\[
\E_{U} \E_0 |L_U(X) - 1|
\le2\PROB\{|U|>a\} + \PROB\{|U|\le a\} \max_{u \in[-a,a]} \E_0
|L_u(X) - 1|
\]
and therefore, using the Cauchy--Schwarz inequality,
\begin{eqnarray*}
1 - \frac{1}{2} \E_{U} \E_0 |L_U(X) - 1|
& \ge&
\PROB\{|U|\le a\}\biggl( 1-\frac{1}{2}\max_{u \in[-a,a]} \E_0
|L_u(X) - 1|\biggr)
\\
& \geq& \PROB\{|U|\le a\}\biggl(1 - \frac12 \max_{u \in[-a,a]}
\sqrt{\E_0 L_u^2(X) - 1}\biggr).
\end{eqnarray*}
Since
\begin{eqnarray*}
L_u(x) & = & \frac{1}{N} \sum_{S \in\cC} \frac{1}{(1-\rho)^{k/2}}
\exp\biggl(- \sum_{i \in S} \frac{(x_i - \sqrt{\rho} u)^2}{2
(1-\rho)} - \sum_{i \notin S} \frac{x_i^2}{2}\biggr) \exp
\Biggl(\sum_{i =1}^n \frac{x_i^2}{2} \Biggr) \\
& = & \frac{1}{N} \sum_{S \in\cC} \frac{1}{(1-\rho)^{k/2}} \exp
\biggl(\sum_{i \in S} \frac{x_i^2}{2} - \frac{(x_i - \sqrt{\rho}
u)^2}{2 (1-\rho)} \biggr),
\end{eqnarray*}
we get
\begin{eqnarray*}
\E_0 L_u^2(X)
& = &\frac{1}{N^2} \sum_{S, S' \in\cC} \frac{1}{(1-\rho)^k} \E_0
\exp\biggl(\sum_{i \in S \cap S'} X_i^2 - \frac{(X_i - \sqrt{\rho}
u)^2}{1-\rho}\\
&&\hspace*{120pt}{} + \sum_{i \in S \Delta S'} \frac{X_i^2}{2} - \frac
{(X_i - \sqrt{\rho} u)^2}{2(1-\rho)}\biggr) \\
& = &\frac{1}{N^2} \sum_{S, S' \in\cC} \frac{1}{(1-\rho)^k (2 \pi
)^{n /2}} \\
&&\hspace*{44pt}{}  \times\int_{-\infty}^{+\infty} \exp\biggl(\sum_{i \in S
\cap S'} \frac{x_i^2}{2} - \frac{(x_i - \sqrt{\rho} u)^2}{1-\rho}
\\
&&\hspace*{104.5pt}{} - \sum_{i \in S \Delta S'} \frac{(x_i - \sqrt{\rho} u)^2}{2(1-\rho
)} - \sum_{i \notin S \cup S'} \frac{x_i^2}{2}\biggr) \,dx.
\end{eqnarray*}
It is easy to check that
\[
\frac{x_i^2}{2} - \frac{(x_i - \sqrt{\rho} u)^2}{1-\rho} = \frac
{\rho u^2}{1+\rho} - \frac{1+\rho}{2(1-\rho)} \biggl(x_i - \frac{2
\sqrt{\rho} u}{1+\rho} \biggr)^2,
\]
which implies
\begin{eqnarray*}
\E_0 L_u^2(X)
& = &\frac{1}{N^2} \sum_{S, S' \in\cC}
\frac{\exp(({\rho u^2}/({1+\rho})) |S \cap S'|)}{(1-\rho)^k (2 \pi)^{n /2}} \\
&&\hspace*{43pt}{}  \times\int_{-\infty}^{+\infty} \exp\biggl(- \sum_{i \in S
\cap S'} \frac{1+\rho}{2(1-\rho)} \biggl(x_i - \frac{2 \sqrt{\rho
} u}{1+\rho} \biggr)^2 \\
&&\hspace*{103.2pt}{} - \sum_{i \in S \Delta S'} \frac{(x_i -
\sqrt{\rho} u)^2}{2(1-\rho)} - \sum_{i \notin S \cup S'} \frac
{x_i^2}{2}\biggr) \,dx \\
& = &\frac{1}{N^2} \sum_{S, S' \in\cC} \frac{\exp(
({\rho u^2}/({1+\rho})) |S \cap S'|)}{(1-\rho)^k} \biggl(\frac
{1-\rho}{1+\rho}\biggr)^{|S \cap S'|/2}\\
&&\hspace*{43pt}{} \times(1-\rho)^{k -
|S \cap S'|} \\
& \leq&\frac{1}{N^2} \sum_{S, S' \in\cC}
\exp\biggl(\biggl(\frac{\rho u^2}{1+\rho}- \frac12\log(1-\rho
^2)\biggr) |S \cap S'| \biggr),
\end{eqnarray*}
which concludes the proof.
\end{pf}

We now apply Theorem~\ref{thmlower} to a~few examples. The theorem
converts the
problem into a~purely combinatorial question and~\cite{combin} offers
various estimates for the moment generating function of $Z$ which we
may use for our purposes.

\subsubsection{Nonoverlapping sets}

Consider first the simplest case when $\cC$ contains $N$ disjoint
sets of size $k$.
%
\begin{cor}
\label{cordisjoint}
Let $\cC$ be the class of all sets of size $k$.
If
\[
\nu(\rho) \le\frac{\log(N)}{k},
\]
then the Bayes risk satisfies $R^*\ge0.3$, and $R^* \to1$ if $\rho
\ll\min(1, \log(N)/k)$ or if $(1-\rho) N^{2/k} \to\infty$.
\end{cor}
\begin{pf}
Clearly, the size $Z$ of the overlap of two randomly chosen elements
of $\cC$ equals zero with probability $1-1/N$ and $k$ with
probability $1/N$. Thus,
\[
\EXP e^{\nu Z}-1 = (1/N)(e^{\nu k} -1) \le(1/N)e^{\nu k},
\]
which is bounded by $1$ if $\nu\le\log(N)/k$. The first part
then follows from the second part of Theorem~\ref{thmlower}. For the
second part, we need to find $a~\to\infty$ such that
$\nu_a~k - \log N \to-\infty$.
(Note that in this case the upper bound above tends to zero.)
First assume that $\rho\ll\min(1, \log(N)/k)$. In that case, $\nu
_a~\sim\rho a^2$, so it suffices to take $a~\to\infty$ slowly enough
that $\rho a^2 \ll\min(1, \log(N)/k)$. Next assume that $b := \log
(1-\rho) + 2 \log(N)/k \to\infty$. In this case, we have $\nu_a~\leq a^2 - (1/2) \log(1-\rho)$, and we simply choose $a~\to\infty$
slowly enough that $a^2 -b/2 \to-\infty$.
\end{pf}

\subsubsection{$k$-intervals}

Consider the class of all $k$-intervals. The situation is similar to
that of nonoverlapping sets.
(In fact, since this class of $k$-intervals contains $[n/k]$
nonoverlapping sets of size $k$,
we could immediately deduce a~lower bound via Corollary~\ref{cordisjoint}.)
%
\begin{cor}
\label{corkint}
Let $\cC$ be the class of all $k$-intervals.
If
\[
\nu(\rho) \le\frac{\log(n/(2k))}{k},
\]
then the Bayes risk satisfies $R^*\ge0.3$, and $R^* \to1$ if $\rho
\ll\min(1, \log(n/k)/k)$ or if $(1-\rho)(n/k)^{2/k} \to\infty$.
\end{cor}
\begin{pf}
For two $k$-intervals chosen independently and uniformly at random,
\[
\mathbf{P}\{|S \cap S'| = \ell\} = \frac2N\qquad\forall\ell= 1,\ldots, k.
\]
%
Thus,
\[
\EXP e^{\nu Z}-1 = \frac2N \Biggl(\sum_{\ell=1}^k e^{\nu\ell}
-k\Biggr) \le\frac{2k}N e^{\nu k},
\]
and proceed as in the proof of Corollary~\ref{cordisjoint}, using the
fact that $N \leq n$.
\end{pf}

\subsubsection{$k$-sets}

Consider the class of all sets of size $k$.
%
\begin{cor}
\label{corksets}
Let $\cC$ be the class of $k$-sets.
If
\[
\frac{k^2}{n} \le\frac{\ln2}{\exp(\nu(\rho)) -1},
\]
then the Bayes risk satisfies $R^*\ge0.3$, and $R^* \to1$ if either
$k^2/n \to\infty$ and $\rho k^2/n \to0$,
or
$(1-\rho)n^2/k^4 \to\infty$.
\end{cor}
%
%
\begin{pf}
By~\cite{combin}, Proposition 3.4, which uses
negative association,
%
\[
\EXP e^{\nu Z} \le\biggl((e^\nu-1) \frac{k}{n}+1\biggr)^k
\le\exp\biggl((e^\nu-1) \frac{k^2}{n}\biggr),
\]
where the last expression is bounded by $2$ under the postulated
condition, and tends to 1
if either $k^2/n \to\infty$ and $\nu k^2/n \to0$, or $k^2/n \to0$
and $e^\nu k^2/n \to0$. First assume that $k^2/n \to\infty$ and
$\rho k^2/n \to0$. By choosing $a~\to\infty$ slowly enough that
$\rho a^2 k^2/n \to0$ we ensure that $\nu_a~k^2/n \to0$. Next assume
that $b := \log(1-\rho) - 2 \log(k^2/n) \to\infty$. Since $\nu_a~\leq a^2 -(1/2) \log(1-\rho)$, it suffices to take $a~\to\infty$
slowly enough that $a^2 -b/2 \to-\infty$ to ensure that $e^\nu k^2/n
\to0$. The result then follows from Theorem
\ref{thmlower}.
\end{pf}

\subsubsection{Perfect matchings}

Consider now the example of perfect matchings described in the
\hyperref[intro]{Introduction}. Here $k=\sqrt{n}$. Once again, Theorem~\ref{thmlower}
applies and implies that testing is impossible for moderate values of
$\rho$.
%
\begin{cor}
\label{cormatch}
Let $\cC$ be the class of all perfect matchings.
If $\rho\le1/2$, the Bayes risk satisfies $R^*\ge0.3$.
Also, $R^* \to1$ if $\rho\to0$.
\end{cor}
\begin{pf}
The random variable $Z$ for this class is considered by
\cite{combin}, who prove that
\[
\EXP e^{\nu Z} \le
\biggl((e^\nu-1) \frac{1}{\sqrt{n}}+1\biggr)^{\sqrt{n}}
\le e^{e^\nu-1}.
\]
This is bounded by $2$ whenever $\nu\le1+\ln\ln2$, which is
satisfied whenever $\rho\le1/2$, and tends to 1 if $\nu\to0$.
We then apply Theorem~\ref{thmlower}.
\end{pf}

\subsubsection{Spanning trees}

A~similar argument applies for the class of all spanning trees of
a~complete graph with $k+1$ vertices [and $n=(k+1)k/2$ edges] as described
in the \hyperref[intro]{Introduction}.
%
\begin{cor}
\label{corsp}
Let $\cC$ be the class of all spanning trees.
If $\rho\le0.4$,
then the Bayes risk satisfies $R^*\ge0.15$.
We also have $R^* \to1$ if $\rho\to0$.
\end{cor}
\begin{pf}
It is shown in~\cite{combin} that
\[
\EXP e^{\nu Z} \le
\biggl((e^\nu-1) \frac{2}{k+1}+1\biggr)^{k}
\le e^{2(e^\nu-1)},
\]
which is bounded by $13/4$ whenever
$\nu\le1+\ln((\ln(13/4))/2)$, which is satisfied whenever $\rho\le
0.4$, and tends to 1 if $\nu\to0$. We then apply Theorem~\ref{thmlower}.
\end{pf}

\section{Some near-optimal tests}
\label{secupper}

We already know that the likelihood ratio test is optimal in the
Bayesian setting.
We study here other tests for multiple reasons. First,
the likelihood ratio test seems difficult to compute in most
situations. Second, the likelihood ratio test is heavily dependent on
the prior we choose---here, the uniform distribution on the class.
The third, and perhaps most important, reason is that
it is difficult to obtain directly upper bounds for the (worst-case) risk
of the likelihood ratio test whereas the tests considered below
are easier to analyze and often yield near-optimal performance.
Whenever we obtain an upper bound for the risk of a~test that
matches the lower
bounds developed in the previous section, we have a~full understanding
of the limitations and possibilities of detection for the
particular case considered, and this is our main goal in this paper.

We
consider the squared-sum test, which corresponds to the ANOVA test in
the detection-of-means setting, the generalized likelihood ratio test
(GLRT) and a~goodness-of-fit (GOF) test, as well as some variants.
We say that a~test is \textit{near-optimal} for a~certain setting if it
achieves the information bound for that setting to first order.

\subsection{The squared-sum test}

One of the simplest tests is based on the observation that
the magnitude of the squared-sum $(\sum_{i=1}^n X_i)^2$ may be
substantially different under the null and alternative hypotheses due
to the higher correlation under the latter.

Indeed, under $\PROB_0$, $(\sum_{i=1}^n X_i)^2$ is distributed as $n
\chi_1^2$, while for
any $S\subset\{1,\ldots,n\}$ with $|S|=k$, under
$\PROB_S$, $(\sum_{i=1}^n X_i)^2$ has the same distribution
as $(n + \rho k(k-1)) \chi_1^2$; in fact, under the more general
correlation model~(\ref{model-general}), this is a~(stochastic) lower
bound. This immediately leads to the following result.
%
\begin{prp}
\label{prpsq}
Let $\cC$ be an arbitrary class of sets of size $k$
and suppose that $\rho k^2/n \to\infty$ in~(\ref{model-general}).
If $t_n$ is such that $t_n\to\infty$ but $t_n = o(\rho k^2/n)$,
then the test which rejects the null hypothesis
if $(\sum_{i=1}^n X_i)^2 > n t_n$
has a~
worst-case
risk converging to zero.
However, any test based on $(\sum_{i=1}^n X_i)^2$ is powerless if
$\rho k^2/n \to0$ in~(\ref{model}).
\end{prp}

In Corollary~\ref{corksets}, we saw that reliable detection of
$k$-sets is
impossible if $k^2/n \to\infty$ and $\rho k^2/n \to0$. Here we see
that, when $\rho k^2/n \to\infty$, the\vadjust{\goodbreak} squared-sum test is
asymptotically powerful. Hence, the following statement:

\begin{quote}
{\normalsize\textit{The squared-sum test is near-optimal for detecting
$k$-sets in the regime where $k^2/n \to\infty$.}}
\end{quote}
On the other hand, in the regime $k^2/n \to0$, the squared-sum test is
powerless even if $\rho= 1$.
The test does not require knowledge of $\rho$, though knowing $\rho$ allows
one to choose the threshold $t_n$ in an optimal fashion; if $\rho$ is
unknown, we simply choose $t_n \to0$ very slowly.

\subsection{The generalized likelihood ratio test}
\label{secglrt}

In this section we investigate the performance of the generalized
likelihood ratio test (GLRT). We show that for parametric classes such
as $k$-intervals, the test is near-optimal. However, for the
nonparametric class of $k$-sets,
the test performs poorly in some regimes.

By definition, the GLRT rejects for large values of $\max_{S \in\cC}
Z_S/\E_0 Z_S$, or simply $\max_{S \in\cC} Z_S$ when all the sets in
the class $\cC$ are of same size, since~$\E_0 Z_S$ only depends on
the size of $S$. Hence, the GLRT is of the form
\[
f(X)=0  \quad\mbox{if and only if}\quad
\max_{S \in\cC} X^T (\bI- \bA_S^{-1}) X \leq t
\]
for some appropriately chosen $t$. We immediately notice that the GLRT
requires knowledge of $\rho$

Our analysis of the GLRT is based on Lemma~\ref{lemqf}, which
provides the
distribution of the quadratic form $X^T (\bI- \bA_S^{-1}) X$ under
the null $\PROB_0$ and under the alternative $\PROB_S$.
Under the null we need to control the maximum of such quadratic forms
over $S \in\cC$, which we do using exponential concentration
inequalities for chi-squared distributions.

\subsubsection{The GLRT for $k$-intervals and other parametric classes}

Recalling Corollary~\ref{corkint}, when detecting $k$-intervals all
tests are
asymptotically powerless when $\rho\ll\min(1, \log(n/k)/k)$. We
assume for concreteness that
$k/\log n \to\infty$, for otherwise detecting $k$-intervals for very
small $k$ has more to do with detecting $k$-sets. We state a~general
result that applies for classes
of small cardinality.
%
\begin{prp} \label{prpglrt-small} Consider a~class $\cC$ of sets of
size $k$, with cardinality $N
\to\infty$ such that $\log(N)/k \to0$. When $\rho k / \log N
\to\infty$, the generalized likelihood ratio test with threshold
value $t = - \rho k + \rho\sqrt{5 k \log N} + 2\log N$ has
worst-case
risk tending to zero.
\end{prp}
\begin{pf}
We first bound the probability of Type I error. Indeed, under the null,
by Lemma~\ref{lemqf} and its proof, we can decompose
\[
X^T (\bI- \bA_S^{-1}) X
= -\frac\rho{1-\rho} C_S + \frac{\rho(k-1)}{1 + \rho(k-1)} D_S,
\]
where $C_S \sim\chi_{k-1}^2$ and $D_S \sim\chi_1^2$. Hence,
\[
\max_{S \in\cC} X^T (\bI- \bA_S^{-1}) X \leq-\rho\min_{S \in
\cC} C_S + \max_{S \in\cC} D_S.
\]
It is well known that the maximum of $N$ standard normals is bounded
by
$\sqrt{2 \log N}$
with probability tending to 1 as $N \to
\infty$. Hence, the second term on the right-hand side is bounded by
$2 \log N$ with high probability. For the first term, we combine the
union bound and Chernoff's bound to obtain, for all $a\le1$,
%
\begin{eqnarray}\label{glrt-small-eq1}
\PROB_0\Bigl\{\min_{S \in\cC} C_S < a~(k-1)\Bigr\}
&\leq& N \mathbf{P}\{\chi_{k-1}^2 < a(k-1)\} \nonumber\\[-8pt]\\[-8pt]
&\leq& N \exp\biggl(-\frac{(k-1)}2 (a~- 1 - \log a)\biggr).
\nonumber
\end{eqnarray}
Using the fact that $a~- 1 - \log a~\sim\frac12 (1-a)^2$ when $a~\to
1$, the right-hand side tends to zero when $a~= 1 - \sqrt{(5/k)\log
N}$. We arrive at the conclusion that the GLRT with threshold $t =
- \rho k + \rho\sqrt{5 k \log N} + 2 \log N$ has probability of Type
I error tending to zero.

Now consider the alternative under $\P_S$. By Lemma~\ref{lemqf} and
Chebyshev's inequality,
\[
X^T (\bI- \bA_S^{-1}) X \geq- \rho k - \rho s_k \sqrt{k} + \rho k /s_k
\]
with high probability when $s_k \to\infty$. We then conclude by the
fact that the right-hand side is larger than $t$ when $s_k \to\infty$
sufficiently slowly.
\end{pf}

Comparing the performance of the GLRT in Proposition~\ref{prpglrt-small} with the
lower bound for $k$-intervals in Corollary~\ref{corkint}, we see that the
GLRT is near-optimal for detecting $k$-intervals. This is actually
the case for all parametric classes we know of.

\subsubsection{The GRLT for $k$-sets and other nonparametric classes}
\label{secglrt-nonparametric}

Consider now the example of the class of all $k$-sets.
Compared to the previous section, the situation here is different in that
$N$, the size of the class $\cC$, is much larger. For example, for
$k$-sets, $N = {n \choose k}$, and
therefore $\log(N)/k \to\infty$ with $n \to\infty$. The
equivalent\vspace*{1pt} of Proposition~\ref{prpglrt-small} for this regime is the
following:
%
\begin{prp} \label{prpglrt-large}
Consider a~class $\cC$ of sets of size $k$, with cardinality $N \to
\infty$ such that $\log(N)/k
\to\infty$. When $\eta:= (1-\rho) N^{2/k} (\log N)/k \to0$, the
generalized likelihood
ratio test with threshold value $t = -(\log N)/\sqrt{\eta}$ has
worst-case
risk tending to zero.
\end{prp}
\begin{pf}
We follow the proof of Proposition~\ref{prpglrt-small}. The only
difference is in~(\ref{glrt-small-eq1}), where we now need $a~\to0$ and that
right-hand side tends to zero when $\log a~+ 2 (\log N)/k \to-\infty
$. Choose $a~= N^{-2/k} \sqrt{\eta}$, obtaining that, with\vadjust{\goodbreak} high probability,
%
\begin{equation} \label{glrt-large-eq1}
\max_{S \in\cC} X^T (\bI- \bA_S^{-1}) X \leq-\frac\rho{1-\rho}
N^{-2/k} k \sqrt{\eta} + 2 \log N.
\end{equation}
As before, with high probability under $\P_S$,
%
\begin{equation} \label{glrt-large-eq2}
X^T (\bI- \bA_S^{-1}) X \geq- \rho k,
\end{equation}
so we only need to check that the threshold $t$ is larger than the
right-hand side in~(\ref{glrt-large-eq1}) and smaller than the
right-hand side in~(\ref{glrt-large-eq2}), which is the case by the
assumptions we made.
\end{pf}

Notice that in Proposition~\ref{prpglrt-large} the condition on $\rho
$ implies
that $\rho\to1$, which is much stronger than what the squared-sum
test requires when $k^2/n \to\infty$.
For $k$-sets, $N = {n \choose k}$---so that $\log N = k \log(n/k) +
O(k)$---and the requirement is that $(1-\rho) (n/k)^2 \log(n/k) \to
0$, which is substantially stronger than what the lower bound obtained
in Corollary~\ref{corksets} requires. Moreover, if we restrict $\rho
$ to be bounded away from $1$, then the GLRT may be powerless.
%
%
\begin{theorem}
\label{thmglrt-bad}
Let $\cC$ be the class of all $k$-sets. If $\rho< 0.6$ and $k=
o(n^{0.7})$, the GLRT has a~
Bayes
risk bounded away from zero.
\end{theorem}

The proof is in the \hyperref[app]{Appendix}.

In view of Theorem~\ref{thmglrt-bad}, the GLRT is clearly suboptimal
when in
the situation stated there, and compares very poorly with the
squared-sum test, which is asymptotically powerful if $\rho k^2/n \to
\infty$ as seen in Proposition~\ref{prpsq}. We do not know of any
other situation
where the GLRT fails so miserably.

\subsection{A~localized squared-sum test}

While the GLRT is near-optimal for detecting objects from a~parametric
class such as $k$-intervals, it needs knowledge of $\rho$.
However, a~simple modification solves this drawback. Indeed, consider
the following ``local'' squared-sum test:
\[
f(X)=0  \quad\mbox{if and only if}\quad
\max_{S \in\cC} \biggl(\sum_{i \in S} X_i\biggr)^2 \leq t
\]
for some appropriate threshold $t$.
%
\begin{prp} \label{prplocal-sum}
Consider a~class $\cC$ of sets of size $k$, with cardinality $N
\to\infty$ such that $\log(N)/k \to0$. When $\rho\gg\log(N)/k$
in~(\ref{model-general}), the local squared-sum test with threshold $t
= 2 k \log N$ has
worst-case
risk tending to zero.
\end{prp}
\begin{pf}
The proof is quite straightforward. Indeed, under the null, for any $S$
of size $k$ we have $\sum_{i \in S} X_i \sim\cN(0, k)$ so that
\[
\max_{S \in\cC} \biggl(\sum_{i \in S} X_i\biggr)^2 \leq t\vadjust{\goodbreak}
\]
with probability tending to 1. Under an alternative~(\ref{model-general}), $S$ denoting the anomalous set of variables, we have
\[
\P\biggl( \biggl(\sum_{i \in S} X_i\biggr)^2 \ge t\biggr) \ge\P\bigl( \bigl(k + k(k-1) \rho
\bigr) \chi_1^2 \ge t\bigr) \to1,
\]
when $\rho\gg\log(N)/k$.
\end{pf}

Specializing this result to the case of $k$-intervals leads to the
following statement (which ignores logarithmic factors):

\begin{quote}
{\normalsize\textit{The localized squared-sum test is near-optimal for
detecting $k$-intervals in the regime where $\log(n)/k \to0$.}}
\end{quote}

\textit{When $k$ is unknown.} We might only know that some interval is
anomalous, without knowing the size of that interval. In that case,
multiple testing at each $k$ using the local squared-sum test yields
adaptivity. Computationally, this may be done effectively by computing
sums in a~multiscale fashion as advocated in~\cite{MGD}. In fact, here
it is enough to compute the sums over all \textit{dyadic}
intervals---since each interval $S$ contains a~dyadic interval of
length at least $|S|/4$---and this can be done in $3 n$ flops in a~recursive fashion.

\subsection{A~goodness-of-fit test}

By now, the parametric case is essentially solved, with the local
squared-sum test being not only near-optimal but also computable in
polynomial time
(in $n$ and $k$)
for the case of $k$-intervals, for example. In the nonparametric case,
so far, the story is not complete. We focus on the class of all
$k$-sets. There we know that the squared-sum test is near-optimal if
$k^2/n \to\infty$. If $k^2/n \to0$, it has no power, and we only
know that the GLRT works when $(1-\rho) (n/k)^2 \log(n/k)\to0$,
which does not match the rate obtained in Corollary~\ref{corksets}.
Worse than that, it is not clear whether computing the GLRT is possible
in time polynomial in $(n, k)$.
We now show that a~simple
goodness-of-fit (GOF)
test performs (almost) as desired.

The basic idea is the following.
Let $H_i = \Phi^{-1}(X_i)$, where $\Phi$ is the standard normal distribution
function. Under the null, the $H_i$'s are i.i.d. uniform in $(0,1)$.
Under an alternative with anomalous set denoted by $S$, the $X_i, i
\in S$ are closer together, especially since we place ourselves in the
regime where $\rho\to1$. More precisely, we have the following.
%
\begin{lem} \label{lemclose}
Suppose $X_1,\ldots, X_k$ are
zero-mean, unit-variance
random variables satisfying $\Cov(X_i, X_j) \geq\rho> 0$, for all $i
\neq j$. Let $\overline{X}$ denote their average. Then for any $t > 0$,
\[
\PROB\bigl\{ \# \{i\dvtx |X_i - \overline{X}| > t\} \geq k/2\bigr\}
\leq\frac{2(1-\rho)}{t^2}.
\]
\end{lem}
\begin{pf}
Let $\Lambda:= \sum_{i \neq j} \Cov(X_i, X_j) \geq k(k-1) \rho$.
Elementary calculations show that
\[
\E\biggl[\frac1k \sum_i (X_i - \overline{X})^2\biggr]
= 1 -\frac1k -\frac{\Lambda}{k^2} \leq(1 -1/k)(1 - \rho) \leq
1-\rho.
\]
By Markov's inequality, we then have
\[
\PROB\biggl\{ \frac1k \sum_i (X_i - \overline{X})^2 > t^2/2
\biggr\}
\leq\frac{2(1-\rho)}{t^2}.
\]
The statement follows from observing that
\[
\# \{i\dvtx |X_i - \overline{X}| > t\} \geq k/2  \quad\Rightarrow\quad
\frac1k \sum_i (X_i - \overline{X})^2 > t^2/2.
\]
\upqed\end{pf}

The idea, therefore, is detecting unusually high concentrations of
$H_i$'s, which is a~form of GOF test for the uniform distribution.
Under a~general correlation model as in~(\ref{model-general}), with
Lemma~\ref{lemclose} we see that the concentration will happen over an
interval of length slightly larger than $\sqrt{1-\rho}$. This is
apparent from Lemma~\ref{lemrepresent} under the simple correlation model
(\ref{model}).

Choose an integer $m$ such that $m \gg(n/k^2) \log(n/k^2)$ and partition
the interval $[0,1]$ into $m$ bins of length $1/m$, denoted $I_s,  s =
1,\ldots, m$. Let $B_s = \# \{i\dvtx H_i \in I_s\}$ be the bin
counts---thus, we are computing a~histogram. Then consider the
following GOF test:
\[
f(X)=0  \quad\mbox{if and only if}\quad
\max_{s = 1,\ldots, m} B_s \leq t,
\]
where $t$ is some threshold.
%
%
\begin{prp} \label{prpgof}
Consider the class $\cC$ of all $k$-sets in the case where $k^2/n \to
0$ and $k/\log n \to\infty$. In the GOF test above, choose $m$ such
that $(n/k^2) \log n \ll m \ll n/\log n$. When $(1-\rho)^{1/2} \ll
1/m$ in~(\ref{model-general}), the resulting test with threshold $t =
n/m + \sqrt{3 n \log(m)/m}$ has
worst-case
risk tending to zero.
\end{prp}
\begin{pf}
Bernstein's inequality, applied to the binomial distribution, gives that
\[
\PROB_0\bigl\{B_s > n/m + b \sqrt{n/m}\bigr\} \leq\exp\bigl[- (b^2/2)/\bigl(1 + (b/3)
\sqrt{m/n}\bigr)\bigr].
\]
This and the union bound imply that, indeed,
\[
\PROB_0\Bigl\{\max_s B_s > t\Bigr\} \to0.
\]

Consider now an alternative of the form~(\ref{model-general}), with
$S$ denoting the anomalous set. Let
\[
I := \{i \in S\dvtx |X_i - \overline{X}_S| \leq1/m\},\qquad  \overline
{X}_S := \frac1k \sum_{i \in S} X_i.\vadjust{\goodbreak}
\]
Though the set $I$ is random, by Lemma~\ref{lemclose} and the fact that
$(1-\rho)^{1/2} \ll1/m$, we have that
\[
\PROB_S\{|I| \geq k/2\} \to1.
\]
Define the event $Q := \{-a~\leq\overline{X}_S \leq a\}$ for some $a~> 0$. Note that, since the variance of $\overline{X}_S$ is bounded by
1, $\PROB(Q^c) \leq2 (1 - \Phi(a))$. Define $\tilde{H}_S = \Phi
^{-1}(\overline{X}_S)$. On~$Q$, using a~simple Taylor expansion, we have
\[
|H_i - \tilde{H}_S| \leq\frac{|X_i - \overline{X}_S|}{\phi(a~+
1/m)} \leq e^{a^2}/m\qquad \forall i \in I,
\]
where $\phi$ denotes the standard normal density function and $a$
is taken sufficiently large. Therefore, when $|I| \geq k/2$ and $Q$
hold, at least $k/2$ of the anomalous $H_i$'s fall in an interval of
length at most $2 e^{a^2}/m$. Since such an interval is covered by at
most $2 e^{a^2}$ bins, by the pigeonhole principle, there is a~bin
that contains $k e^{-a^2}/4$ anomalous $H_i$'s. By Bernstein's
inequality, the same bin will also contain at least $(n-k)/m -\sqrt{3
n \log(m)/m}$ nonanomalous $H_i$'s (with\vspace*{1pt} high probability), so in
total this bin will contain $n/m -k/m -\sqrt{3 n \log(m)/m} + k
e^{-a^2}/4$ points. By our choice of $m$, $k \gg\sqrt{n \log(m)/m}$,
so it suffices to choose $a~\to\infty$ slowly enough that $k e^{-a^2}
\gg\sqrt{n \log(m)/m}$ still. Then, with high probability, there is
a~bin with more than $t$ points.
\end{pf}

Ignoring logarithmic factors, we are now able to state the following:

\begin{quote}
{\normalsize\textit{The GOF test is near-optimal for detecting $k$-sets in
the regime where $k^2/n \to0$ and $k/\log n \to\infty$.}}
\end{quote}

When $k/\log n \to0$, things are somewhat different. There, the GOF
test requires that $(1-\rho) n^{2 k/(k-1)} \to0$, which is still
close to optimal when $k \to\infty$, but far from optimal when $k$ is
bounded (e.g., when $k=2$, the exponent is 4 instead of~2).
Indeed, when $k/\log n \to0$, $m$ needs to be chosen larger than $n$,
and Bernstein's inequality is not accurate. Instead, we use the simple bound
\[
\P\bigl(\operatorname{Bin}(n,p) \ge\ell\bigr) \le2 \frac{ (n p)^\ell}{\ell!}
\qquad\mbox{when } n p \le1/2.
\]
Note that Bennett's inequality would also do. (The analysis also
requires some refinement showing that, with probability tending to 1
under the alternative, one cell contains at least $k$ points.)
Note that in the remaining case, $k = O(1)$, the GLRT is optimal up to
a~logarithmic factor, since it only requires that $(1-\rho) n^2 \log n
\to0$, as seen in Section~\ref{secglrt-nonparametric}. We do not
know whether
a~comparable performance can be achieved by a~test that does not have
access to $\rho$.

\textit{When $k$ is unknown.} In essence, we are trying to detect an
interval with a~higher mean in a~Poisson count setting.\vadjust{\goodbreak} As before, it
is enough to look at dyadic intervals of all sizes, which can be done
efficiently as explained earlier, following the multiscale ideas
in~\cite{MGD}.

\begin{appendix}\label{app}
\section*{\texorpdfstring{Appendix: Proof of Theorem \lowercase{\protect\ref{thmglrt-bad}}}{Appendix: Proof of Theorem 3.1}}

The proof is divided into three steps. The first step formalizes the
fact that we want to prove that (under $H_1$),
the contaminated set has no influence
(with high probability) on the GLRT statistic. The second step
exhibits a~useful high probability event. Finally, in the third step we
show that on this high probability event, the contaminated set has no
influence on the GLRT.

It can easily be seen that for every $S$ of size $k$,
\[
X^T (\bI- \bA_S^{-1}) X
= \frac{\rho}{(1+\rho(k-1))(1-\rho)}
\biggl(\sum_{i,j \in S, i \neq j} X_i X_j - \rho(k-1) \sum_{i \in S}
X_i^2 \biggr).
\]
Introduce the function $g\dvtx\R^k\to\R$ defined by
\[
g(u) = \sum_{i \neq j} u_i u_j - \rho(k-1) \sum_{i} u_i^2
= \Biggl(\sum_{i=1}^n u_i\Biggr)^2 - \bigl(1 + \rho(k-1)\bigr)
\sum_{i=1}^n u_i^2
\]
for $u=(u_1,\ldots,u_k)\in\R^k$.
Denoting, for $x\in\R^n$ and $S\subset\{1,\ldots,n\}$,
the vector of components of $x$ belonging to $S$ by $x|_S$,
we may write the GLRT as
\[
f(x)=0  \quad\mbox{if and only if}\quad  \max_{S \in\cC} g(x |_S)
< t.
\]
Note that by the symmetry of $\cC$ and the test,
\begin{eqnarray*}
R(f) & = & \P_0\Bigl\{\max_{S \in\cC} g(X |_S) \geq t\Bigr\} +
\frac1N \sum_{S' \subset\cC} \P_{S'}\Bigl\{\max_{S \in\cC} g(X
|_S) < t\Bigr\} \\
& = & \P_0\Bigl\{\max_{S \in\cC} g(X |_S) \geq t\Bigr\}
+ \P_{{\{1,\ldots,k\}}}\Bigl\{\max_{S \in\cC} g(X |_S) < t
\Bigr\}.
\end{eqnarray*}
Given $X \sim\cN(0, \bI)$, define the coupling $X'$ as follows: $X_i
= X_i'$ for $i \notin\{1,\ldots,k\}$,
and $X_i, X_i'$ are independent for $i \in\{1,\ldots,k\}$. Note that
$X' \sim\cN(0,\bA_{\{1,\ldots,k\}})$. Then,
no matter what the threshold $t$ is, we have
\begin{eqnarray*}
R(f) & = & \P\Bigl\{\max_{S \in\cC} g(X |_S) \geq t\Bigr\} + \P
\Bigl\{\max_{S \in\cC} g(X' |_S) < t\Bigr\} \\
& \geq& \P\Bigl\{\max_{S \in\cC} g(X |_S) \geq\max_{S \in\cC}
g(X' |_S)\Bigr\}.
\end{eqnarray*}
In the following we show that, with probability tending to $1$, we
have
\[
\max_{S \in\cC} g(X |_S) = \max_{S \in\cC} g(X' |_S),
\]
which then implies that the GLRT is asymptotically powerless.

By Lemma~\ref{lemrepresent}, there exist $U, U_1,\ldots, U_k$ independent
standard normal such that for all $i \in\{1,\ldots,k\}$,
\[
X_i' = \sqrt{\rho}   U + \sqrt{1-\rho}   U_i .
\]
Using the fact that $\max_{i=1,\ldots, k} |U_i| \leq\sqrt{2 \log k}$
with high probability, with
probability tending to 1, we have
\[
X_1',\ldots,X_k' \in[-\zeta, \zeta],
\]
where $\zeta:= \sqrt{2 (1-\rho) \log(\omega_k k)}$ and $\omega_k$
is any sequence such that $\omega_k \to\infty$.

Fix $\gamma> 1$ to be determined later and define $p = \mathbf{P}\{
\zeta\leq U \leq\gamma\zeta\}$ where $U \sim\cN(0,1)$. By the
fact that $X_1,\ldots, X_n$ are i.i.d. standard normal, $Z := \# \{i\dvtx
\zeta\leq X_i
\leq\gamma\zeta\} \sim\operatorname{Bin}(n, p)$, so that $\mathbf{P}\{Z
\geq k\} \to
1$ if $k = o(np)$. When $\gamma$ is bounded away from $1$, this is the
case if $\sqrt{\log k}   k^{2 - \rho} = o(n)$.

In conclusion, we proved that the event
\begin{eqnarray*}
\Omega& = & \bigl\{ X_1',\ldots,X_k' \in(-\zeta, \zeta
),   \exists\alpha_1,\ldots, \alpha_k, \beta_1,\ldots,
\beta_k \in\{1,\ldots,n\}
\mbox{ distinct:}  \\
& &\hspace*{118pt}  X_{\alpha_1},\ldots, X_{\alpha_k}, -X_{\beta_1},\ldots
, -X_{\beta_k} \in( \zeta, \gamma\zeta) \bigr\}
\end{eqnarray*}
has a~probability that tends to $1$ if $\sqrt{\log k}   k^{2 - \rho}
= o(n)$ as long as $\gamma$ is bounded away from 1.


We specify $\gamma= 1/\sqrt{\rho+ (\frac1{k-1} + \rho
)^2}$. Note that, as required, $\gamma$ exceeds and is bounded away
from 1.
Assume that we are on the event $\Omega$. First note that
%
\setcounter{equation}{0}
\begin{eqnarray}\label{eqminordered}
g(X_{\alpha_1},\ldots, X_{\alpha_k}) & \geq& k (k-1)
\zeta^2 - \rho(k-1) k \gamma^2 \zeta^2 \nonumber\\[-8pt]\\[-8pt]
& = & k (k-1) \zeta^2 (1 - \rho\gamma^2 ), \nonumber
\end{eqnarray}
and the same holds for $g(X_{\beta_1},\ldots,X_{\beta_k})$.

Let $S \in\cC$ be such that $S \cap\{1,\ldots,k\} \neq\varnothing
$. We want to show that there exists $S'$ such that $g(X |_{S'}) \geq
g(X' |_S)$. This entails that $\max_{S \in\cC} g(X |_S) \geq\max
_{S \in\cC} g(X' |_S)$, since for $S \cap\{1,\ldots,k\} =
\varnothing$ we have $g(X |_S)=g(X' |_S)$. First remark that we can
assume that
%
\begin{equation} \label{eqsumsquared}
\biggl(\sum_{i \in S} X_i' \biggr)^2 \geq\zeta(k-1) \sqrt{1 - \rho
\gamma^2},
\end{equation}
since otherwise by~(\ref{eqminordered}) we can simply take
$S'=\{\alpha_1,\ldots,\alpha_k\}$. To simplify notation, we may assume
that $1 \in S \cap\{1,\ldots,k\}$. By definition of $\Omega$ and the
fact that $S$ contains at least one index in $\{1,\ldots, k\}$, there
exist $u, v \in
\{1,\ldots, k\}$ such that $X_{\alpha_u}$ and $X_{\beta_v}$ do not
appear in $X'|_S$. We want to show that by replacing $X_1'$ by either
$X_{\alpha_u}$ or $X_{\beta_v}$, in $X'|_S$, one increases the value
of $g$. More precisely, we want to show that
\[
\max\bigl(g\bigl(X_{\alpha_u},
X' |_{S \setminus\{1\}}\bigr), g\bigl(X_{\beta_v}, X' |_{S \setminus\{1\}}\bigr)
\bigr) \geq g(X'|_S).
\]
Then by induction one can show the existence
of the $S'$ described above.\vadjust{\goodbreak}

Note that, for $x \in\R^k$ and $y \in\R$,
\begin{eqnarray*}
&&g(x_1,\ldots,x_{j-1}, y, x_{j+1}, \ldots, x_k) - g(x)
\\
&&\qquad = 2 (y- x_j) \sum_{i \neq j} x_i - \rho(k-1) (y^2 - x_j^2) \\
&&\qquad = (y - x_j) \Biggl(2 \sum_{i=1}^k x_i - \bigl(2 + \rho(k-1) \bigr) x_j -
\rho(k-1) y \Biggr) .
\end{eqnarray*}
Consider the case where $\sum_{i \in S} X_i' > 0$ (the case
$\sum_{i \in S} X_i' < 0$ can be dealt with similarly). Since
$X_{\alpha_u} \geq X_1'$, it suffices to show that
$2 \sum_{i \in S} X_i' \geq(2 + \rho(k-1) ) X_1' + \rho(k-1)
X_{\alpha_u}$,
which follows from
\begin{eqnarray*}
\bigl(2 + \rho(k-1) \bigr) X_1' + \rho(k-1) X_{\alpha_u}
& \leq& (k-1) \zeta
\gamma\biggl( \frac{2}{k-1} + 2 \rho\biggr) \\
& = & 2 (k-1) \zeta\sqrt{1 - \rho\gamma^2} \\
& \leq& 2 \sum_{i \in S} X_i .
\end{eqnarray*}
This concludes the proof.
\end{appendix}

\section*{Acknowledgments}

We thank Omiros Papaspiliopoulos for his illuminating remarks and the
anonymous referees for challenging us to obtain stronger results in the
sparse setting and for pointing out a~mistake in Proposition~\ref{prpgof}.



\printaddresses

\end{document}